\documentclass{amsart}
\usepackage{amssymb}
\usepackage{amsmath}
\usepackage{amsthm}
\usepackage{amscd}
\usepackage{srcltx}
\long\def\comment#1\endcomment{}

\usepackage{color}

\def\mcgb{\MCG (S)_b^\omega}
\def\BBF{\mathcal{BBF}}

\def\bu{\mathbf{U}}

\def\bv{\mathbf{V}}

\def\by{\mathbf{Y}}

\def\MCG{\mathcal{MCG}}

\def\ss{\mathcal{U}}
\def\upss{\Pi \ss /\omega}
\def\upcs{\Pi C_k/\omega}

\def\AM{\mathcal{AM}}

\def\stab{{\rm stab}}
\def\CC{{\mathcal C}}

\newcommand{\Cutp}{{\mathrm{Cutp}}\, }

\newcommand{\tdist}{\widetilde{\mathrm{dist}}}

\newcommand{\dist}{{\mathrm{dist}}}
\newcommand{\dcs}{\dist_{\CC (S)}}

\newcommand {\q}{\mathfrak q} 
\newcommand {\g}{\mathfrak g} 
\newcommand {\pgot}{\mathfrak p}

\newcommand{\la}{\langle}
\newcommand{\ra}{\rangle}

\newcommand{\fh}{\mathfrak{h}}
\newcommand{\fk}{\mathfrak{k}}

\newcommand{\uas}{{$\omega$-almost surely}}
\newcommand{\uass}{$\omega$-almost surely }

\newcommand {\R}{\mathbb{R}} 

\newcommand {\free}{\mathbb{F}} 

\def\calc{\mathcal{C}}   

\def\Z{{\mathbb Z}}

\catcode`\@=11
\long\def\@savemarbox#1#2{\global\setbox#1\vtop{\hsize\marginparwidth
  \@parboxrestore\tiny\raggedright #2}}

\newcommand{\tsh}[1]{\left\{\kern-.9ex\left\{#1\right\}\kern-.9ex\right\}}

\newcommand {\me}{\medskip}

\newtheorem{theorem}{Theorem}
\newtheorem{proposition}[theorem]{Proposition}
\newtheorem{corollary}[theorem]{Corollary}

\newtheorem{lemma}[theorem]{Lemma}

\theoremstyle{definition}
\newtheorem{defn}[theorem]{Definition}
\theoremstyle{definition}

\theoremstyle{remark}
\newtheorem{remark}[theorem]{Remark}

\def\square{\hfill${\vcenter{\vbox{\hrule height.4pt \hbox{\vrule
width.4pt
height7pt \kern7pt \vrule width.4pt} \hrule height.4pt}}}$}

\def\<{\langle}
\def\>{\rangle}

\newcommand\seq[1]{\mathchoice{\mbox{\boldmath$#1$}}{\mbox{\boldmath$#1$}}{\mbox{\boldmath$\scriptstyle#1$}}{\mbox{\boldmath$\scriptstyle#1$}}}

\newcommand\seqrep[1]{\lm_\omega\left(#1_{n}\right)}

\newcommand\ultra[1]{( #1_{n})^{\omega}}

\newcommand{\lm}{{\lim}}

\newcommand{\lio}[1]{\lm_\omega\left(#1\right)}

\newcommand{\co}{\colon\thinspace}

\begin{document}
\title[Homomorphisms into mapping class groups]{Homomorphisms into mapping class groups. An addendum}
\author{Jason Behrstock}\thanks{The research of the first author was
supported in part by PSC-CUNY grant \# 60051-39 40.}
\address{Lehman College,
City University of New York, 
U.S.A.}
\email{jason.behrstock@lehman.cuny.edu}
\author{Cornelia Dru\c{t}u}\thanks{The research of the second author was supported in part by
the ANR project ``Groupe de recherche de G\'eom\'etrie et Probabilit\'es dans
les Groupes''.}
\address{Mathematical Institute,
24-29 St Giles,
Oxford OX1 3LB,
United Kingdom.}
\email{drutu@maths.ox.ac.uk}
\author{Mark Sapir}\thanks{The research of the third author was supported in part by NSF grant DMS-0700811.}
\address{Department of Mathematics,
Vanderbilt University,
Nashville, TN 37240, U.S.A.}
\email{m.sapir@vanderbilt.edu}

\begin{abstract}This is an addendum to \cite{BDS}. We show, using our methods and an auxiliary result of Bestvina-Bromberg-Fujiwara, that a finitely generated group with infinitely many pairwise non-conjugate homomorphisms to a mapping class group virtually acts non-trivially on an $\R$-tree, and, if it is finitely presented, it virtually acts non-trivially on a simplicial tree.
\end{abstract}

\maketitle

The goal of this addendum to \cite{BDS} is to show that our methods
together with a result of Bestvina, Bromberg and Fujiwara
\cite[Proposition~5.9]{BBF}
yield a proof of the following theorem.

\begin{theorem}\label{t1.1} If a finitely presented group $\Gamma$ has infinitely many
pairwise non-conjugate homomorphisms into $\MCG(S)$, then $\Gamma$
virtually splits (virtually acts non-trivially on a simplicial
tree).
\end{theorem}

$\aleph_0$
This theorem is a particular case of a result announced
by D. Groves.\footnote{Groves first announced a version of this
result at MSRI in 2007, \cite{Groves}.
More recent announcements by Groves have included stronger
versions of Theorem~\ref{t1.1}.} From private emails received by the authors,
it is clear that the methods used by Groves are significantly different.
Note that the same new methods allow us to give another
proof of the finiteness of the set of homomorphisms from a group with property~(T) to a mapping class group  \cite[Theorem~1.2]{BDS} which is considerably shorter than our original proof; see
Corollary~\ref{c7} below and the discussion following it. Theorem 1.2 in \cite{BDS} may equally be obtained from Theorem \ref{t1.1} above and the fact that every group with property (T) is a quotient of a finitely presented group with property (T) \cite[Theorem p. 5]{shalom:invent2000}.

The property of the mapping class groups given by  Theorem~\ref{t1.1} can be viewed as another ``rank 1" feature of these groups. In contrast, note that  a recent
result of \cite{LRT} shows that the rank 2 lattice $SL_{3}(\Z)$ contains infinitely many pairwise
non-conjugate copies of the triangle group $\Delta(3, 3, 4) =\la a, b \mid a^3 = b^3 = (ab)^4 = 1\ra$. Also, as was pointed out to us by Kassabov, although the group $SL_3(\Z[x])$ has property (T) \cite{Sh2}, it has infinitely many pairwise non-conjugate homomorphisms into $SL_3(\Z)$ induced by ring homomorphisms $\Z[x]\to \Z$.

\smallskip

The following proposition
contains one of the main auxiliary results in \cite{BBF}
and the key ingredient missed in our treatment of groups with many
homomorphisms into mapping class groups in \cite{BDS}.

\begin{proposition}[Bestvina, Bromberg, Fujiwara, Proposition~5.9 of \cite{BBF}]\label{th3.1}
There exists an explicitly defined finite index torsion-free subgroup $\BBF(S)$ of $\MCG(S)$ such that
the set of all subsurfaces of $S$ can be partitioned
into a finite number of subsets $C_1, C_2,...,C_s$, each of which is an orbit
of $\BBF(S)$, and any two subsurfaces in the same subset
overlap and have the same complexity.
\end{proposition}

The proof of this important result, explained to us by Bestvina, is
surprisingly simple: the subgroup $\BBF(S)$ is the subgroup of mapping
classes from $\MCG(S)$ acting as identity on the factor $\pi_1(S)/B$
over certain characteristic subgroup $B$ of $\pi_1(S)$ of finite index
which is explicitly constructed.

We consider the set of colors $K=\{ 1, 2,...,s\}$, and we color each subsurface of $S$ contained in the subset $C_{i}$ by $i$. Note that the whole surface $S$ has a color which is different from that of
any proper subsurface.

Recall that for every sequence  of subsurfaces $\bu$ from $\upss$
we defined an $\R$-tree $T_\bu$ (see \cite[Notation 4.4]{BDS}) and that there is an
equivariant bi-Lipschitz embedding $\psi$ of $\AM$ into
$\prod_{\bu\in \upss} T_\bu$ (see \cite[Corollary 4.17]{BDS}). Let $C_k$ be the set of all
subsurfaces of $S$ with the given color $k\in K$. Let $\pi_k$ be the
projection of $\prod_{\bu\in \upss} T_\bu$ onto $\prod_{\bu\in
\upcs} T_\bu$.

\begin{remark}\label{rem1.1} By \cite[Lemma 2.1]{BDS}, we have that $\prod_{\bu\in \upss} T_\bu$ can be
written as $$\prod_{k\in K}\prod_{\bu\in \upcs} T_\bu.$$
\end{remark}

In what follows we use the notion of tree-graded space introduced in \cite{DrutuSapir:TreeGraded}.

\begin{theorem}\label{th4.1}
 Consider an arbitrary color $k\in K$ and the image $T_k = \pi_k\psi (\AM )$.

 For every subsurface $\bu \in \upcs$ consider the tree $T'_\bu= T_\bu\times \prod_{\bv\in \upcs\setminus\{\bu\}} \{a_\bv \}$ where $a_\bv $ is the point in $T_\bv$ which is the projection of $\partial\bu$ to $T_\bv$.

The space $T_k$ is tree-graded with respect to $T_\bu'\, $ and with transversal trees reduced to singletons. In particular it is an $\R$-tree.
\end{theorem}

\proof \emph{Step 1.} \quad We prove by induction on $n$ that for any finite subset
$F\subset \upcs$ of cardinality $n$ the projection $\pi_F (\AM )$ of $\AM$ onto the
finite product $\prod_{\bu \in F} T_{\bu }$ is an $\R$-tree. The case $n=1$ is obvious, the case $n=2$ follows from
\cite[Theorem 4.21, (2)]{BDS}, since the subsurfaces in $F$ pairwise overlap. Assume that the statement is proved for $n$ and consider $F\subset \upcs$, $F$ of cardinality $n+1\, $.

Both $\psi (\AM )$ and its projections are geodesic
spaces. For $\psi (\AM )$ this follows from Proposition 4.18, while for projections it follows from the fact that the distance is $\ell^1$. To prove that $\pi_F(\psi (\AM ))$ is a real tree it suffices therefore to prove that it is $0$-hyperbolic, i.e. for every
geodesic triangle its three edges have a common point. By Lemma 4.30 the subset $\pi_F(\psi (\AM ))$ is median, thus it suffices to prove that for an arbitrary triple of points $\seq\nu , \seq\rho , \seq \sigma$ in $\pi_F(\psi (\AM ))$ and every geodesic $\g$ joining $\seq\nu , \seq\rho$ in $\pi_F(\psi (\AM ))$ the median point $\seq\mu\, $ of the triple is on $\g \, $.

Assume that there exist $\bu , \bv$ such that the projection of $\seq\mu $ on $T_\bu \times T_\bv$ is not $(v,u)\, $. Assume that it is $(x,u)$, with $x\neq v$ (the other case is similar).

Consider the projection on the product $\prod_{\by \in F\setminus \{ \bv \}} T_{\by }$. By the inductive hypothesis, $\pi_F(\psi (\AM ))$  projects onto a real tree, in particular there exists $\seq\mu'$ on $\g$ such that its projection on $\prod_{\by \in F\setminus \{ \bv \}} T_{\by }$ coincides with that of $\seq\mu \, $. In particular $\pi_\bu (\mu')= \pi_\bu (\mu) =x\, .$ This implies that the projection on $T_\bu \times T_\bv$ of both $\seq\mu'$ and $\seq\mu$ is $(x,v)\, $ (the unique point with first coordinate $x$). This implies that all coordinates of $\seq\mu'$ and $\seq\mu$ are equal, thus the two points coincide.

Assume now that for every pair $\bu , \bv$ in $F$ the projection of $\seq\mu $ on $T_\bu \times T_\bv$ is $(v,u)\, $. Fix such a pair.
  By the inductive hypothesis and an argument as above there exists $\seq\mu_1\in \g$ such that its projection on $\prod_{\by \in F\setminus \{ \bu \}} T_{\by }$ coincides with that of $\seq\mu$. Similarly there exists $\seq\mu_2\in \g$ such that its projection on $\prod_{\by \in F\setminus \{ \bv \}} T_{\by }$ coincides with that of $\seq\mu$. Then on $T_\bu \times T_\bv$ the point $\seq\mu_1$ projects onto some $(x,u)$ and  $\seq\mu_2$ projects onto some $(v,y)$. This implies that there exists some $\seq\mu'$ on $\g$ between $\seq\mu_1$ and  $\seq\mu_2$ projecting on $T_\bu \times T_\bv$ in $(v,u)\, $. Note that for every $\by \in F\setminus \{ \bu , \bv \}$ the projection of $\seq\mu'$ coincides with that of $\seq\mu_1$ and  $\seq\mu_2$, hence with that of $\seq\mu\, $. It follows that $\seq\mu'= \seq\mu$.

\me

We now prove by induction on $n$ that for any finite $F\subset \upcs$ of cardinality $n$ the projection $\pi_F (\AM )$ of $\AM$ onto the
finite product $\prod_{\bu \in F} T_{\bu }$ is tree-graded with respect to the trees $T_\bu^F = T_\bu \times \prod_{\bv\in F\setminus\{\bu\}} \{a_\bv \}$, where $a_\bv $ is the projection of $\partial\bu$ to $T_\bv\, $. It only remains to prove that $\pi_F (\AM )$ is complete and that it is covered by $T_\bu^F\, .$  Both statements are proved simultaneously when proving that $\pi_F (\AM )$ equals the union $\bigcup_{\bu \in F} T_\bu^F\, .$ Clearly the union is contained in $\pi_F (\AM )\, $. Conversely, consider a point  $x=(x_1,...,x_{n+1})$ in $\pi_F (\AM )\, $. The inductive hypothesis applied to $(x_1,...,x_n)$ and $(x_2,...,x_{n+1})$ implies that for each $n$-tuple there exists $\bu \in F$ such that for every $\bv \neq \bu $ the corresponding coordinate is $\pi_\bv (\bu )$, that is the point in $T_\bv$ which is the projection of $\partial\bu$ to $T_\bv$. Assume that in $(x_1,...,x_n)$ the surface $\bu$ corresponds to the first coordinate, and that in $(x_2,...,x_{n+1})$
 the surface $\bu'$ corresponds to the last coordinate. The projection $(x_1, x_{n+1})$ of $x$ on $T_\bu \times T_{\bu'}$ is either of the form $(\pi_\bu (\bu' ), x_{n+1})$ or of the form $(x_1, \pi_{\bu'} (\bu))\, $. In the first case $x$ is in $\prod_{\bv\in F\setminus\{\bu'\}} \{\pi_\bv (\bu' ) \} \times T_{\bu'}$, in the second $x$ is in $T_\bu \times \prod_{\bv\in F\setminus\{\bu\}} \{\pi_\bv (\bu ) \}\, .$

\me

 \emph{Step 2.} \quad We now prove the statements on $T_k\, $. First we prove that $T_k$ is a real tree, using  an approximation argument similar to that in the proof
that $\AM$ is a median space (\cite[Theorem 4.25]{BDS}). Since $T_k$ is a complete geodesic space, it suffices to prove that it is zero hyperbolic. Thus it suffices to prove that
 for every triple  $\seq\alpha\, ,\seq\beta,\, \seq\gamma$ and $\seq\mu$ its median point, $\seq\mu$ is on any geodesic $\g$
 joining $\seq\alpha$ and $\seq\beta$ in $\psi (\AM)$.

 Assume that the distance from $\seq\mu$ to $\g$ is $\varepsilon>0\, $. Take a finite set of surfaces $F$ s.t.
 the projections of $\seq\alpha\, ,\seq\beta,\, \seq\mu$ in $\prod_{\bu\not\in F} T_\bu$ compose a set of diameter $\varepsilon/4\, $.
 Since the projection on the cartesian product $\prod_{\bu\in F} T_\bu$ is a tree,
 the projection of $\g$ contains that of $\seq\mu$, hence there exists $\seq\mu'$ on
 $\g$ with the same projection as $\seq\mu$ in $\prod_{\bu\in F} T_\bu$.

 Then the distance from $\seq\mu'$ to $\seq\mu$ is $$\sum_{\bu\not\in F}\widetilde{\dist}_\bu (\seq\mu'\, ,\, \seq\mu) \leq \sum_{\bu\not\in F}\left[\widetilde{\dist}_\bu (\seq\mu'\, ,\, \seq\alpha) + \widetilde{\dist}_\bu (\seq\alpha  \, ,\,  \seq\mu)\right] \leq  \varepsilon/4 + \varepsilon/4 = \varepsilon/2 \, .
 $$

 The tree $T_k$ is complete. Consider two points $\seq\mu$ and $\seq\nu$ in $\AM\, .$ There exists a countable family $\calc \subset \upcs$ equal to the set of subsurfaces $\{ \bu \; ;\; \dist_\bu (\seq\mu \, ,\, \seq\nu ) >0 \, \}$. Let $\mathfrak h$ be a hierarchical path joining $\seq\mu$ and $\seq\nu\, $. Let $\seq\alpha$ and $\seq\beta$ be the endpoints of a minimal sub-arc ${\mathfrak h}_\bu$ on $\mathfrak h$ such that $\tdist_\bu (\seq\alpha ,\seq\beta ) = \tdist_\bu (\seq\mu ,\seq\nu )\, .$ Assume that there exists $\bv \neq \bu\, ,\, \bv \in \calc\, ,$ such that $\tdist_\bv (\seq\alpha ,\seq\beta )>0\, .$ Then by projecting $\mathfrak h$ on $T_\bu \times T_\bv$ and using the tree-graded structure of the projection of $\AM$ we obtain that the arc ${\mathfrak h}_\bu$ has a strict sub-arc of endpoints $\seq\alpha'$ and $\seq\beta'$ such that $\tdist_\bu (\seq\alpha' ,\seq\beta' ) = \tdist_\bu (\seq\mu ,\seq\nu )\, .$ This contradicts the minimality of ${\mathfrak h}_\bu
 \, .$ It follows that for every  $\bv \neq \bu\, ,\, \bv \in \calc\, ,$ $\tdist_\bv (\seq\alpha ,\seq\beta )= 0\, .$ Hence ${\mathfrak h}_\bu$ is entirely contained in a factor $T_\bu \times \prod_{\bv \in \calc\, ,\, \bv \neq \bu } \{ a_\bv\} \, .$ Since given any subsurface $\bv \neq \bu$ the arc ${\mathfrak h}_\bu$ contains points with first coordinate distinct from the projection of $\bv$ on $T_\bu$ it follows that $a_\bv$ is the projection of $\bu$ on $T_\bv \,$. Hence ${\mathfrak h}_\bu$ is contained in the tree $T_\bu'$, and the arcs ${\mathfrak h}_\bu$ with $\bu \in \calc$ cover ${\mathfrak h}$ up to a subset of zero measure.\endproof

As an immediate consequence of Theorem \ref{th4.1}, we obtain the following, which also immediately follows from the main result of \cite{BBF}.

\begin{corollary}\label{c6} There exists an equivariant embedding of $\AM$ into a
finite product of $\R$-trees.
\end{corollary}

Now let $\Gamma$ have infinitely many pairwise non-conjugate
homomorphisms into $\MCG(S)$. Theorem \ref{th4.1} and Proposition \ref{th3.1}
imply that $\Gamma$ has a finite index subgroup $\Gamma'$ that acts
on the $\R$-trees $T_k$ for each $k\in K$, further, since the global
action is non-trivial (i.e., without a global fixed point) at least one of the actions on a factor tree is non-trivial.

Corollary \ref{c6} and the standard argument of Bestvina and Paulin \cite{Bes1, Pau} imply the following

\begin{corollary} \label{c7}
If a finitely generated group $\Lambda$ has infinitely many pairwise non-conjugate homomorphisms into the group
$\MCG(S)$, then $\Lambda$ has a subgroup of index at most $|K|$ which
is not an $F\R$-group (i.e., acts non-trivially on an $\R$-tree).
\end{corollary}

Since in a group with property (T) every subgroup of finite index has
property $F\R$ \cite{Pau}, Corollary 6.3 in \cite{BDS} follows from Corollary \ref{c7}.

It is still unknown if every finitely generated group acting
non-trivially on an $\R$-tree also acts non-trivially on a simplicial
tree.  In order to obtain such an action in our case, we apply the theorem of Bestvina and Feighn below.

\begin{defn} \label{d5}
Given an action of a group on an $\R$-tree, an arc $\g_0$ is called {\em stable} if the stabilizer of every non-trivial subarc of $\g_0$ is the same as the stabilizer of $\g_0$.

The action is called {\em stable} if every arc $\g$ contains a non-trivial stable subarc $\g_0$.
\end{defn}

\begin{theorem}{\rm (Bestvina-Feighn, \cite[Theorem
    9.5]{BestvinaFeighn:stableactions}).}\label{bf} Let $G$ be a finitely presented group with a nontrivial, minimal, and stable action on an $\R$-tree $T$. Then either
(1) $G$ splits over an extension $E$-by-cyclic subgroup where $E$ is the stabilizer of a non-trivial arc of $T$, or
(2) $T$ is a line. In the second case, $G$ has a subgroup of index at most 2
  that is the extension of the kernel of the
action by a finitely generated free abelian group.
\end{theorem}

In order to show stability of the action, as in \cite{DrutuSapir:splitting}, we describe stabilizers of pairs of points and of tripods in
$T_k$.

The following Lemmas \ref{stab2} and \ref{stab3} describing stabilizers of arcs and tripods have similar proofs as Lemmas 5.14 and 5.15 in the main text.

\begin{lemma}\label{stab2} There exists a constant $N=N(S)$, such that
    if $\seq\mu$ and $\seq\nu$ are distinct points in $\AM$ that are
    not in the same piece, then the  stabilizer $\stab (\seq\mu , \seq\nu )$
    is the extension of a finite subgroup of cardinality at most $N$ by an
    abelian group.
\end{lemma}

\proof By hypothesis, for every representatives $(\mu_n)$ and
$(\nu_n)$ of  $\seq\mu$ and $\seq\nu$
respectively, the following is satisfied:
\begin{equation}\label{infdcs}
\lim_\omega \dcs(\mu_n ,\nu_n)=\infty\, .
\end{equation}

Let $\seq g=\ultra g$ be an element in $\stab (\seq\mu , \seq\nu )$.
Then
$$\delta_n(\seq g) =\max (\dist (\mu_n , g_n \mu_n)\, ,\, \dist (\nu_n
, g_n \nu_n))$$ satisfies $\delta_n(\seq g) = o(d_n)$. Let $\q_n$ be a
hierarchy path joining $\mu_n$ and $\nu_n$ and let $\bar{\mu}_n$ and
$\bar{\nu}_n$ be points on $\q_n$ at distance $\varepsilon d_n$ from
$\mu_n$ and respectively $\nu_n$. By hypothesis for $\varepsilon$
small enough $\lim_\omega \dcs (\bar{\mu}_n , \bar{\nu}_n) =\infty$.
Thus there exist $\tilde{\mu}_n$ and $\tilde{\nu}_n$ on $\q_n$
between $\bar{\mu}_n$ and $\bar{\nu}_n$ and at respective
$\CC(S)$-distance 3 from them. Denote by $\q_n'$ the sub-arc of
$\q_n$ between $\tilde{\mu}_n$ and $\tilde{\nu}_n$.

Divide
$\q_n'$ into three consecutive sub-arcs that shadow geodesics
in $\CC (S)$ of equal length $\frac{\dcs
(\tilde{\mu}_n,\tilde{\nu}_n)}{3}$. Let us show that there exists a point $x=(x_n)^\omega$ on the first part and a point $y=(y_n)^\omega$  on the third part which are at distance $O(1)$ from $g\pgot'$ (the points do not depend on $g$).

All large domains on $\q_n'$ are \uass large domains for $g_n \q_n$. Suppose
that the whole surface $S$ is the only large domain of a part $\pgot_n$ of $\q_n'$ of
size $O(d_n)$. Then we can take a projection of $\g_n$ and $g_n\g_n$
to the curve complex $\CC(S)$ and deduce from the hyperbolicity of
$\CC(S)$ that the geodesics $\pgot_n'$ and $g_n\pgot_n'$ are at
$\CC(S)$-distance $O(1)$ $\omega$-a.s. Thus in that case we can take points $(x_n)^\omega$ and $(y_n)^\omega$ arbitrarily.

Suppose that such a large domain in $\pgot_n$
cannot be found $\omega$-a.s. Note that the distance between the
entry points of $\g_n'$ and $g_n\g'_n$ into large domains $S'\subset
S$ are at $\CC(S)$-distance 1. Thus in this case we can take $(x_n)^\omega$ and $(y_n)^\omega$ to be the entrance points of the geodesic into large domains.

Obviously $\lim_\omega \dcs (x_n, y_n)=\infty$.

For every $\seq g=\ultra g\in\stab (\seq\mu , \seq\nu)$ we define a
sequence of translation numbers. Since $x_n$ is \uass at distance
$O(1)$
from a point $x_n'$ on $g_n q_n$, define $\ell_x (g_n)$ as
$(-1)^\epsilon \dcs (x_n , g_n x_n) $, where $\epsilon =0$ if $x_n'$
is nearer to $g_n \mu_n$ than $g_n x_n $ and $\epsilon =1$
otherwise.

Let $\ell_x \co \stab (\mu , \nu ) \to \Pi \R /\omega $ defined by
$\ell_x (\seq g)= \left(\ell_x (g_n) \right)^\omega$. It is easy to
see
that $\ell_x$ is a quasi-morphism,
that is
\begin{equation}\label{qm}
\left|  \ell_x (\seq g\seq h)- \ell_x (\seq g)-\ell_x
(\seq h) \right| \leq _{\omega} O(1)\, .
\end{equation}

It follows that $\left|  \ell_x \left( [\seq g,\seq h] \right)
\right| \leq _{\omega} O(1) \, .$

The above and a similar argument for $y_n$ imply that for every
commutator, $\seq c = \seqrep c$, in the stabilizer of $\seq \mu$ and
$\seq\nu$,
$\dcs (x_n , c_n x_n)$ and $\dcs (y_n , c_n y_n)$ are at most $O(1)$.
 Bowditch's
acylindricity result \cite[Theorem 1.3]{Bowditch:tightgeod} and Lemma~2.1 
imply that the set of
commutators of $\stab (\seq\mu, \seq\nu)$ has uniformly bounded
cardinality, say, $N$. Then any finitely generated subgroup $G$ of
$\stab (\seq\mu, \seq\nu)$ has conjugacy classes of cardinality at
most
$N$, i.e. $G$ is an $FC$-group \cite{Neumann:finiteconj}. By
\cite{Neumann:finiteconj}, the set of all torsion elements of $G$ is
finite, and the derived subgroup of $G$ is finite of cardinality
$\le N(S)$ (by Lemma~2.13). 
\endproof


\begin{lemma}[Lemma 2.20, (2), in \cite{DrutuSapir:splitting}]
\label{inv} Let $\free$ be a tree-graded space.  For every non-trivial
geodesic $\g$ in the tree obtained by collapsing
non-trivial pieces, $T=\free /\!  \!  \approx$, there
exists a non-trivial geodesic $\pgot$ in $\free$ such that its
projection on $T$ is $\g$, and such that given an isometry $\phi$ of
$\free$ permuting the pieces, the isometry $\tilde \phi$ of $T$
induced by $\phi$ fixes $\g$ pointwise if and only if $\phi$ fixes the
set of cutpoints $\Cutp(\pgot)$ pointwise.
\end{lemma}

The quotient tree $\AM /\! \!  \approx$ is described in
\cite[Lemma 3.8]{BDS}.

\begin{lemma}\label{stab3}
Let $\widetilde{\seq\mu}_1$, $\widetilde{\seq\mu}_2$ and
$\widetilde{\seq\mu}_3$ be three points in the quotient tree $\AM /\!
\!  \approx$ which form a non-trivial tripod.
Then the stabilizer $\stab (\widetilde{\seq\mu}_1 , \widetilde{\seq\mu}_2,
\widetilde{\seq\mu}_3 )$ in $\mcgb$ is a finite subgroup of
cardinality at most $N=N(S)$.

\end{lemma}

\proof
For every $i\in \{1,2,3\}$ let $\g_i$ denote the geodesic joining
$\widetilde{\seq\mu}_j$ and $\widetilde{\seq\mu}_k$ in $T$, where
$\{i,j,k\}= \{1,2,3\}$, and let $\pgot_i$ denote
a geodesic in $\AM$ associated to  $\g_i$ by Lemma \ref{inv}.  By eventually replacing the endpoints of
$\pgot_i$ with cut-points in their interiors we may assume that the
three geodesics $\pgot_1\, ,\, \pgot_2$ and $\pgot_3$ compose a
triangle in $\AM$ of vertices ${\seq\alpha}$, ${\seq\beta}$ and
${\seq\gamma}$.  Note that the elements in $\stab (\widetilde{\seq\mu}_1,
\widetilde{\seq\mu}_2, \widetilde{\seq\mu}_3)$ fix
point-wise all the cut-points of all the geodesics $\pgot_i\,$.
Since the set of cut-points does not change, we may replace the three
geodesics by three paths $\fh^{i}$, each of which is an
ultralimit of a
sequences of hierarchy paths, $\lio{\fh_n^i}$, with the property
that the endpoints of $\fh_n^1,
\fh_n^2, \fh_n^3$ are in the set of vertices of a triangle, $\mu_n^1,
\mu_n^2, \mu_n^3$.  Each $\fh_n^i$ projects onto a geodesic
$\gamma_n^i$ in the curve complex $\CC(S)$, and according to \cite[Lemma
4.15]{BDS} we also have $\lio {\mathrm{length} (\gamma_n^i)} = \infty \,$.

By hyperbolicity of $\CC(S)$ for every $a>0$ there exists $b>0$ such that for any triple of points $x,y,z\in\CC(S)$ the intersection of the three
$a$--tubular neighborhoods of geodesics $[x,y]$, $[y,z]$, and
$[z,x]$ is a set $C_{a}(x,y,z)$ of diameter at most $b$. In
particular for every $n$ the three $a$--tubular neighborhoods of the
geodesics $\gamma_n^1, \, \gamma_n^2,\, \gamma_n^3$ intersect in a set
$C_n$ of diameter at most $b$.  Fix an $\epsilon>0$ and
consider a sub-path $\fk_n^1$ of $\fh_n^1$ such that the limit path $\fk_1 =\lio{\fk_n^1}$ has endpoints at $\tdist$-distance $\epsilon$ and $2\epsilon$ from $\seq\mu_2$. Consider a (sufficiently) large proper domain $Y_n^2$ for $\fk_n^1$. If no proper large
domain exists for $\fk_n^1$ (i.e. the only large domain for this hierarchy path is $S$) then
pick instead a marking $\rho_{n}^{1}$ on $\fk_n^1$. Since
we started with a non-trivial tripod, for
$\epsilon$ small enough the sub-arc $\fk_1$ is at positive
$\widetilde{\dist}$-distance from $\fh_2$, hence $Y_n^2$ is \uass not
a large domain of $\fh_n^2$ (or, in the second case,
$\rho_{n}^{1}$ is not at uniformly bounded $\CC(S)$-distance from  $\fh_n^2$).
Therefore $Y_n^2$  is a large domain of $\fh_n^3$
(respectively $\rho_{n}^{1}$ is at $\CC(S)$-distance $O(1)$ from
$\fh_n^3$).  Let $g =
(g_n)^\omega$ be an element of $\stab (\widetilde{\seq\mu}_1,
\widetilde{\seq\mu}_2, \widetilde{\seq\mu}_3)$. Consider any geodesic
quadrangle with two of the opposite edges being
$\fh_n^1$ and $g_n\fh_n^1$. Since $\fk_1$ is at positive $\widetilde{\dist}$-distance both from $\seq\mu_2$ and from $\seq\mu_3$, the
domain $Y_n^2$ (or the marking $\rho_{n}^{1}$) cannot be at uniformly bounded $C(S)$-distance from the edges $[\mu_n^2 \, ,\, g_n \mu_n^2]$ and $[\mu_n^3 \, ,\, g_n \mu_n^3]$ of the quadrangle. Thus,  $Y_n^2$ can only be a large domain of $\fh_n^1$ and
$g_n\fh_n^1$
(respectively, only these two edges contain points
at $\CC(S)$-distance $O(1)$ from $\rho_{n}^{1}$).
A similar argument shows that $Y_n^2$ is a large domain of (or  $\rho_{n}^{1}$ is at $\CC(S)$-distance $O(1)$ from) $g_n\fh_n^3\, $.

In a similar manner we take a sub-path $\fk_n^2$ of $\fh_n^2$ such that the limit path $\fk_2 =\lio{\fk_n^2}$ has endpoints at $\widetilde{\dist}$-distance $\epsilon$ and $2\epsilon$ from $\seq\mu_1$, we fix $Y_n^1$ proper large domain for $\fk_n^2$ (or a marking $\rho_n^1$ on $\fk_n^2$ if no such domain exists). Then we show that $Y_n^1$ is also a large domain for $\fh_n^3$, $g_n\fh_n^2$ and $g_n\fh_n^3$ (respectively $\rho_n^1$ is at $\CC(S)$-distance $O(1)$ from these paths).
Likewise, we find a large domain $Y_n^3$ for $\fh_n^1$ and
$\fh_n^2$ and their translations by $g_n$  (or a marking $\rho_n^3$ at $C(S)$-distance $O(1)$ from all these paths).

Let $\widehat{\fh}_n^1$ be the sub-arc of $\fh_n^1$ between the
sub-arcs corresponding to the domains $Y_n^2$ and $Y_n^3$ (respectively the sub-arc between the
markings $\rho_n^2$ and $\rho_n^3$), and $\widehat{\gamma}_n^1$ its
projection into the complex of curves. Note that $\widehat{\gamma}_n^1$ is a
sub-arc of $\gamma_n^1$. Likewise consider $\widehat{\fh}_n^i$ and $\widehat{\gamma}_n^i$ for $i=2,3$. The
set $C_n$ equals also the intersection of the three $a$-tubular
neighborhoods of the geodesics $\widehat{\gamma}_n^1, \,
\widehat{\gamma}_n^2,\, \widehat{\gamma}_n^3\, $.  Indeed, it clearly
contains this intersection.  On the other hand, the existence of a
point in $C_n$ not in the intersection would imply, up to reindexing, the
existence of a point in $\gamma_n^1 \setminus \widehat{\gamma}_n^1$ at
finite $\CC(S)$-distance from both $\gamma_n^2$ and $\gamma_n^3\, .$
All elements in $\gamma_n^1 \setminus \widehat{\gamma}_n^1$ are
projections in $\CC(S)$ of sub-arcs of $\fh_n^1$ with limits at
$\tdist$-distance at most $2\epsilon$ from either $\seq\mu_2$ or
$\seq\mu_3\, $.  For $\epsilon$ small enough these limits are
therefore at positive $\tdist$-distance from either $\fh_2$ or
$\fh_3$, hence the ultralimit of the
$\CC(S)$-distance of the corresponding sequence of
sub-arcs of $\fh^1_n$ either to $\fh_n^2$ or to $\fh_n^3$ is $\infty$.

 The translation $g_n C_n$ is the intersection of the three $a$-tubular neighborhoods of the geodesics $g_n\widehat{\gamma}_n^1, \, g_n\widehat{\gamma}_n^2,\, g_n\widehat{\gamma}_n^3\, $. For every $i\, $, on the path $g_n\fh_n^i$ the two large domains $Y_n^j$ and $g_n Y_n^j$ occur such that the corresponding sub-arcs have limits at $\tdist$ zero. Then with an argument as above it can be proved that $g_n C_n$ is also the intersection of three $a$-tubular neighborhoods of geodesics of $\CC(S)$ joining the projections of $Y_n^1, Y_n^2,Y_n^3\, $. It follows that $C_n$ and $g_n C_n$ are at Hausdorff distance at most $D=D(S)\, $.

According to the above, there exists $\lambda_{n}$ satisfying $\lio
{\lambda_n}=\infty$ and points $\alpha_n$ on $\gamma_n^1$ at
distance at least $2\lambda_n$ from the projections of the domains
$Y_n^2, Y_n^3, g_nY_n^2, g_nY_n^3$ and
 such that $g_n \alpha_n$ is at distance $O(1)$ from
$\alpha_n\, $.  We pick $\beta_n$ on $\gamma_n^1$ at distance
$\lambda_n$ from $\alpha_n \, $.  Then $g_n \beta_n$ is on $g_n
\gamma_n^1$ at distance $\lambda_n$ from $g_n\alpha_n\, $.

Since $\beta_n$ is on a geodesic between $\alpha_n$ and the projection
of $Y_n^2$, say, and both endpoints are at distance $O(1)$ from $g_n
\gamma_n^1$ it follows that there exists $\beta_n'$ on $g_n \gamma_n$
at distance $O(1)$ from $\beta_n$.  It follows that $\beta_n'$ is at
distance $\lambda_n +O(1)$ from $g_n \alpha_n$, hence it is at
distance $O(1)$ from $g_n \beta_n\, $.  We have thus obtained
$\alpha_n$ and $\beta_n$ at $\CC(S)$-distance $\lambda_n$ such that
$g_n\alpha_n$ is at $\CC(S)$-distance $O(1)$ from $\alpha_n$, and $g_n\beta_n$ is at $\CC(S)$-distance $O(1)$ from $\beta_n$.  It now follows from
Bowditch's acylindricity result \cite[Theorem
1.3]{Bowditch:tightgeod} and \cite[Lemma 2.1]{BDS} that
$\stab (\widetilde{\seq\mu}_1,
\widetilde{\seq\mu}_2, \widetilde{\seq\mu}_3)$
has uniformly bounded
cardinality.\endproof

\begin{lemma}\label{l2} Let $\BBF(S)_b^\omega$ be the subset in $\Pi_b \MCG(S)/\omega$ composed of elements $(x_i)^\omega$ with $x_i\in
\BBF(S)$ \uas.  Then
$\BBF(S)_b^\omega$ is a torsion-free subgroup of index $|\MCG(S)/\BBF(S)|$ in $\Pi_b \MCG(S)/\omega$. Moreover $\BBF(S)_b^\omega$ acts on
each $T_k$ faithfully.
\end{lemma}

\proof Only the last statement requires a proof.  An element $g_\omega
= (g_n)^\omega$ in $\Pi_b \MCG(S)/\omega$ which acts by fixing
$T_{k}$ pointwise  must fix pointwise $T_{\bu}$ for each $\bu \in \Pi
C_k/\omega$.  In particular, for each $U\in C_k$ the mapping class
$g_n$ fixes \uass its boundary $\partial U$.  Each $C_k$ contains a
pair of subsurfaces whose boundaries fill the surface, and the only
mapping classes which fix a pair of filling curves are those of finite
order (uniformly bounded by the complexity of $S$).  Hence only finite
order elements of $\Pi_b \MCG(S)/\omega$ can be in the kernel of the
homomorphism $\BBF(S)_b^\omega \to \mathrm{Isom} (T_k)$.  Since
$\BBF(S)_b^\omega$ is torsion free, the proof is complete.\endproof

\begin{corollary} \label{c3}
A finitely generated $F\R$ group $\Lambda$ cannot have infinitely many pairwise non-conjugate homomorphisms into the group $\BBF(S)$.
\end{corollary}

Let $\bu=(U_i)^\omega$ be an element of $\upcs$ and let $T_\bu'$ be the corresponding sub-tree in $T_k\, .$

\begin{lemma}\label{l5}
\begin{enumerate}
  \item The stabilizer in $\BBF(S)^\omega_b$ of a non-trivial arc in
  $T_\bu'$ has a homomorphism onto a
  (finite of cardinality at most $N=N(S)$)-by-abelian subgroup $A$ of
   $\Pi_b \MCG(U_i)/\omega$. The kernel $W$ of that homomorphism acts identically on $T_\bu'$.
      \me

  \item The stabilizer in  $\BBF(S)^\omega_b$ of a non-trivial tripod in $T_\bu'$ has a homomorphism onto a finite of cardinality at most $N=N(S)$ subgroup of   $\Pi_b \MCG(U_i)/\omega$; the kernel of that homomorphism is $W$.
\end{enumerate}
\end{lemma}

\proof Let $g$ be an element in $\BBF(S)^\omega_b$ stabilizing a
non-trivial arc $\mathfrak h$ in $T_\bu'$. Then $g$ stabilizes $\bu$.  Indeed, we have $gT'_\bu=T'_{g\bu}$.  If
$g\bu\ne \bu$ then $T'_\bu$ and $T'_{g\bu}$ intersect in more than
one point (since they both contain $\mathfrak h$),
which is impossible since these trees are the pieces in a tree-graded
structure.  Therefore the stabilizer of $\mathfrak h$ in
 $\BBF(S)^\omega_b$ must stabilize $\bu$. Hence there exists a homomorphism from that stabilizer to $\MCG_b (\bu)$ whose kernel fixes $T_\bu'$ pointwise. By Lemma \ref{stab2} the image $A$ of that homomorphism is (finite of cardinality at most $N=N(S)$)-by-abelian.

If instead of the stabilizer of an arc in $T_\bu'$ we consider the
stabilizer of a tripod, the argument is similar, except that we use
Lemma \ref{stab3} instead of \ref{stab2}.\endproof

\begin{lemma} \label{l5.1} Let $\Lambda$ be a finitely generated group
with infinitely many pairwise non-conjugate homomorphisms into
$\MCG(S)$.  Then $\Lambda$ contains a subgroup
$\Lambda'$ of index at most $|K|$ which acts on each of  the limit trees $T_k$.
Moreover, each of the actions of $\Lambda'$ on $T_k$ is stable.
\end{lemma}

\proof
That $\Lambda$ contains a subgroup $\Lambda'$ of index at
most $|K|$ which acts on each of the trees $T_k$ follows immediately
from Corollary \ref{c6}.  We now prove that these actions are stable.

By Theorem \ref{th4.1}, the tree $T_k$ is a tree-graded space with pieces the trees
$T'_\bu$ and with all the
transversal trees consisting of singletons. Hence every geodesic $\g$ in $T_k$ is covered, up to a subset of measure zero, by (countably many) non-trivial arcs in trees $T'_\bu$.

Consider an arbitrary $\bu \in \upss$ and the intersection of
$\Lambda'$ with the stabilizer of $T_\bu'$ in $\MCG_b^\omega(S)$, denoted by $\Lambda_\bu$. In view of Lemma \ref{l5}, in order to prove stability it suffices to prove that
stabilizers in $\Lambda_\bu$ of non-trivial arcs in $T_\bu'$ satisfy the
ascending chain condition.  Consider the homomorphism $\pi \colon
\Lambda_\bu\to \Pi_b \MCG(U_i)/\omega$ defined in Lemma \ref{l5}.  The
stabilizer in $\Lambda_\bu$ of a non-trivial arc $\fh$ in $T_\bu'$ is
the inverse image by $\pi$ of the stabilizer of $\fh$ in $\pi
(\Lambda_\bu)$.  Thus it is enough to prove that stabilizers of arcs
in $\pi (\Lambda_\bu )$ satisfy the ascending chain condition.
According to Lemma \ref{l5} the stabilizers of arcs in
$\pi(\Lambda_\bu )$ are (finite of cardinality at most
$N(S)$)-by-abelian, and stabilizers of tripods are finite of
cardinality at most $N(S)\, $.  According to
\cite[Lemma 2.35]{DrutuSapir:splitting} an arc with stabilizer in $\pi
(\Lambda_\bu )$ of order larger than $(N+1)!$ is stable.  (Note that
the hypothesis in Lemma 2.35 that the group acting be finitely
generated is not needed in the proof.)  The ascending chain condition
is obviously satisfied on the set of stabilizers of sub-arcs of order
at most $(N+1)!  \, $.\endproof

Now Theorem \ref{t1.1} follows from Theorem \ref{bf} and Lemma \ref{l5.1}.


\def\cprime{$'$}
\providecommand{\bysame}{\leavevmode\hbox to3em{\hrulefill}\thinspace}
\providecommand{\MR}{\relax\ifhmode\unskip\space\fi MR }
\providecommand{\MRhref}[2]{%
  \href{http://www.ams.org/mathscinet-getitem?mr=#1}{#2}
}
\providecommand{\href}[2]{#2}

\end{document}